\documentclass[12pt,a4paper]{article}
\usepackage{amsmath,amssymb,amsfonts}

\title{\bf Periodic continued fractions and Kronecker symbols}

\author{Kurt Girstmair}
\date{}

\makeatletter
\let\@@maketitle=\maketitle
\def\maketitle{\def\thispagestyle##1{\relax}\@@maketitle}
\makeatother
\newtheorem{theorem}{Theorem}
\newtheorem{prop}{Proposition}
\newtheorem{lemma}{Lemma}

\def\BE{\begin{equation}}
\def\EE{\end{equation}}
\def\BD{\begin{displaymath}}
\def\ED{\end{displaymath}}
\def\BA{\begin{array}}
\def\EA{\end{array}}
\def\BEA{\begin{eqnarray}}
\def\EEA{\end{eqnarray}}
\def\BI{\bibitem}

\def\phi{\varphi}
\def\EPS{\varepsilon}

\def\MOD{\: \mbox{mod} \:}

\def\MB{\mbox}
\def\LD{\ldots}
\def\OV{\overline}

\def\WT{\widetilde}

\def\MN{\medskip\noindent}
\def\STOP{\hfill$\blacksquare$}

\def\JS#1#2{ \left( \frac{#1}{#2} \right) }

\begin{document}
\maketitle

\begin{abstract}
\noindent
We study the Kronecker symbol $\JS st$ for the sequence of the convergents $s/t$ of a purely periodic
continued fraction expansion. Whereas the corresponding sequence of Jacobi symbols is always periodic,
it turns out that the sequence of Kronecker symbols may be aperiodic. Our main result describes the
period length in the periodic case in terms of the period length of the sequence of Jacobi symbols
and gives a necessary and sufficient condition for the occurrence of the aperiodic case.
\end{abstract}

\section*{1. Introduction and main results}

Let $[a_0, a_1, a_2, \LD]$ be the regular continued fraction expansion
of an irrational number $z$. We consider the {\em convergents} $s_k/t_k$, $k\ge 0$, of this
expansion. They are defined by the well-known recursion formulas
\BE
\label{1.0}
s_{-1}=1,\ t_{-1}=0,\ s_0=a_0,\ t_0=1,
\EE
and
\BE
\label{1.1}
  s_k=a_ks_{k-1}+s_{k-2}, \enspace t_k=a_kt_{k-1}+t_{k-2},
\EE
for $k\ge 1$. Then
\BD
  s_k/t_k=[a_0,\LD,a_k],\ k\ge 0,
\ED
is the $k$th convergent of $z$ (see \cite{Hu},
p. 250). Note that $t_k$ is a positive integer for $k\ge 0$.

In the papers \cite{Gi1} and \cite{Gi2} we investigated the {\em Jacobi symbol} $\JS{s_k}{t_k}$
in the periodic case, i.e., for a quadratic irrational $z$. Since this symbol is defined
only for odd denominators $t_k$, we defined $\JS{s_k}{t_k}=*$ if $t_k$ is even. It turned out that the sequence of
Jacobi symbols $\JS{s_k}{t_k}$, $k\ge 0$, is periodic with a period length $L=dl$, where $l$
is the smallest possible period length of $[a_0, a_1, a_2, \LD]$ and $d$ is a divisor of $8$ or $12$.
We called this sequence the {\em Jacobi sequence} of $z$.

The natural generalization of the Jacobi symbol $\JS st$ for arbitrary co-prime integers $s, t$, $t\ge 1$, is the {\em Kronecker symbol}.
It coincides with the Jacobi symbol if $t$ is odd. If $t=2^jt'$, where $j\ge 1$ and $t'$ is an odd natural number, one defines
\BD
  \JS st=\JS s2^j\JS s{t'},
\ED
with
\BD
  \JS s2=\left\{
           \begin{array}{rl}
             1, & \MB{ if } s\equiv\pm 1 \mod 8; \\
            -1, & \MB{ if } s\equiv\pm 3 \mod 8.
           \end{array}
         \right.
\ED
(see \cite{Co}, p. 28 ff.). The Kronecker symbol shares many properties with the Jacobi symbol, for instance, the reciprocity law
\BE
\label{1.2}
 \JS st=\EPS(s',t')\JS ts,
\EE
where $s$ and $t$ are co-prime, $s=2^js'$, $t=2^lt'$ with odd natural numbers $s'$, $t'$, and
\BD
\EPS(s',t')=\left\{
              \begin{array}{rl}
                -1, & \MB{ if } s',t' \MB{ are both } \equiv 3 \mod 4; \\
                1, & \MB{ otherwise}.
              \end{array}
            \right.
\ED
So one might think that the periodicity of the Jacobi sequence can be generalized to the corresponding sequence
of Kronecker symbols, which we call the {\em Kronecker sequence} of $z$. This, however, is not the case, as our main result shows.

In this paper we restrict ourselves to the {\em purely} periodic case since the mixed periodic one seems to be much more difficult.
So let $z=[\OV{a_0,\LD,a_{l-1}}]$, where $l$ has been chosen smallest possible. In \cite{Gi1} we have seen that the Jacobi sequence
of $z$ is purely periodic with a period length $L=dl$, $d\ge 1$, such that
\BE
\label{1.4}
  D_L=\left(
        \begin{array}{cc}
         s_{L-1} & s_{L-2} \\
         t_{L-1} & t_{L-2} \\
        \end{array}
      \right)
   \equiv I \mod 4.
\EE
Here $I$ is the $2\times 2$ unit matrix and the congruence has to be understood entry-by-entry. Further, we may assume that $L$ is even
(in fact, $d$ can be chosen as a divisor of $8$ or $12$, but we do not require this in what follows, since our periods are not always shortest possible).

In this setting suppose that $D_L=I+2^m U$, $m\ge 2$, where not all entries
of
\BE
\label{1.6}
  U=\left(
      \begin{array}{cc}
        x & y \\
        u &v \\
      \end{array}
    \right)
\EE
are even. Suppose, further, that $u=2^e u'$, with $e\ge 0$ and $u'$ odd. A convergent
$s_k/t_k$ of $z$ is called {\em critical} with respect to $L$, if $k\le L-1$, $s_k\equiv 3 \mod 4$ and $t_k\equiv 0\mod 2^{m+e}$
(since $m\ge 2$, the last-mentioned  condition requires $t_k\equiv 0\mod 4$, of course.) Now our main result reads as follows.

\begin{theorem} 
\label{t1}
Let the above notations hold, in particular, let $z$ be purely periodic and $L$ be a period length of the Jacobi sequence with the above properties.
Suppose, further, that no critical convergent with respect to $L$ exists. Then the Kronecker sequence $\JS{s_k}{t_k}$, $k\ge 0$,
is purely periodic with period length $L$ or $2L$. If there is, conversely, a critical convergent with respect to $L$, then the Kronecker
sequence is aperiodic.

\end{theorem} 

\noindent
{\em Remark.} In Proposition \ref{p2} we describe the cases of period length $L$ and $2L$ of the theorem precisely.
As a rule, one finds more examples with periodic Kronecker sequences than with aperiodic ones.

\MN
{\em Examples.}  1. For $z=[\OV{1,2,3}]=(4+\sqrt{37})/7$ we may choose $L=6$, $m=2$. Here $u=21$, so $e=0$. There is no critical convergent among
$s_0/t_0$, \LD, $s_5/t_5$, but the convergent $s_1/t_1=3/2$ has the effect that the Kronecker sequence has only period length $2L=12$.

2. In the case $z=[\OV{1,2,5}]=(7+\sqrt{82})/11$ we may choose $L=12$, $m=2$. Here $u$ is odd, so $e=0$. Hence the convergent $s_7/t_7=975/608$ is critical with respect to $L$
and the Kronecker sequence is aperiodic.

3. An aperiodic example with $e>0$ is $z=(5+\sqrt{85})/10=[\OV{1,2,2}]$, where $L=36$ works with $m=3$ and $e=1$. Therefore, a critical convergent $s_k/t_k$ must satisfy $t_k\equiv 0\mod 16$. The convergent $s_6/t_6=91/64$ has this property.

\section*{2. The reciprocal Jacobi sequence}
An obvious way to generalize the results concerning the Jacobi sequence consists in the generalization of the auxiliary results needed for this purpose. It turns out,
however, that this is impossible, as the following example shows. If $s/t=[a_0,\LD,a_k]$ is a rational number, then the Jacobi symbol $\JS st$ (which equals $*$, if $t$ is even)
depends only on the residue classes of $a_0,\LD,a_k$ mod 4. No result of this kind can hold for the Kronecker symbol.
Indeed, let $s/t=3/2^j$ and $s'/t'=3/(7\cdot 2^j)$, where $j\ge 3$ is odd. We have
\BD
  s/t=[0,(2^j-2)/3,1,2] \MB{ and } s'/t'= [0,(7\cdot 2^j-2)/3, 1,2].
\ED
Here $(2^j-2)/3\equiv (7\cdot 2^j-2)/3 \mod 2^{j+1}$. The Kronecker symbol, however, takes the values $\JS st=-1$ and $\JS{s'}{t'}=1$.

Hence our approach to the Kronecker sequence differs from the above strategy of generalizing auxiliary Jacobi results. Instead, we consider the {\em reciprocal} Jacobi sequence
$\JS{t_k}{s_k}$, $k\ge 0$ (with $\JS{t_k}{s_k}=*$ if $s_k$ is even). Then we use the reciprocity law (\ref{1.2}) in order to obtain the values of the Kronecker symbol in the case where $t_k$ is even.
For this purpose we need the following proposition.

\begin{prop} 
\label{p1}
Let $z=[\OV{a_0,\LD,a_{l-1}}]$.
Suppose that $L=dL$, $d\ge 1$, is an even period length of the Jacobi sequence of $z$
such that {\rm (\ref{1.4})} holds. Then the reciprocal Jacobi sequence $\JS{t_k}{s_k}$, $k\ge 0$, is purely periodic with the same period length $L$.

\end{prop} 

\noindent
{\em Proof.} From the identity
\BE
\label{2.0}
  \left(
    \begin{array}{cc}
      s_{k+L} & s_{k+L-1} \\
      t_{k+L} & t_{k+L-1} \\
    \end{array}
  \right)
  =D_L\cdot
  \left(
    \begin{array}{cc}
      s_{k} & s_{k-1} \\
      t_{k} & t_{k-1} \\
    \end{array}
  \right),\enspace k\ge 0,
\EE
and (\ref{1.4}) we obtain
\BE
\label{2.2}
s_{j+L}\equiv s_j \mod 4,\enspace t_{j+L}\equiv t_j \mod 4 \MB{ for all } j\ge -1.
\EE
Hence $\JS{t_{k+L}}{s_{k+L}}=\JS{t_k}{s_k}=*$ if $s_k$ is even, $k\ge 0$.

If both $s_k$ and $t_k$ are odd, we have
\BE
\label{2.4}
\JS{t_{k+L}}{s_{k+L}}=\EPS(t_{k+L},s_{k+L})\JS{s_{k+L}}{t_{k+L}},
\EE
by quadratic reciprocity. Now (\ref{2.2}) shows $\EPS(t_{k+L},s_{k+L})=\EPS(t_k,s_k)$.
Since $L$ is a period length of the Jacobi sequence, $\JS{s_{k+L}}{t_{k+L}}=\JS{s_{k}}{t_{k}}$.
Accordingly, (\ref{2.4}) says
\BD
\JS{t_{k+L}}{s_{k+L}}=\EPS(t_{k},s_{k})\JS{s_{k}}{t_{k}}.
\ED
Now quadratic reciprocity shows $\JS{t_{k+L}}{s_{k+L}}=\JS{t_{k}}{s_{k}}$.

There remains the case $t_k$ even, $s_k$ odd. Then $k\ge 1$, since $t_0=1$.  We use the notation of the proof of Theorem 5 in \cite{Gi1} and put
$s=s_{k-1}$, $t=t_{k-1}$, $p=a_k$, $q=1$, $m=s_k$, $N=t_k$ and $\delta=(-1)^{k-1}$. Since $N$ is even, $t$ must be odd. Two cases have to be distinguished:

{\em Case 1}: $s$ is odd.
\\ Then the said theorem yields
\BD
  \JS{\delta t}{s}\JS pq\JS{-\delta N}{m}=\EPS(s,q,m)
\ED
with $\EPS(s,q,m)=-1$, if two of the numbers $s,q,m$ are $\equiv 3\mod 4$, and $\EPS(s,q,m)=1$, otherwise. This gives
\BD
  \JS Nm = \JS{-\delta}{m}\JS{\delta}s \EPS(s,1,m)\JS ts.
\ED
Now $\EPS(s,1,m)=\EPS(s,m)$, and quadratic reciprocity implies
\BD
 \JS Nm = \JS{-{\delta}}{m}\JS{\delta}s \EPS(s,m)\EPS(t,s)\JS st,
\ED
i.e.,
\BE
\label{2.6}
\JS {t_k}{s_k} = \JS{-\delta}{s_k}\JS{\delta}{s_{k-1}} \EPS(s_{k-1},s_k)\EPS(t_{k-1},s_{k-1})\JS{s_{k-1}}{t_{k-1}}.
\EE
In the same way we obtain
\BD
\JS {t_{k+L}}{s_{k+L}} = \JS{-\delta'}{s_{k+L}}\JS{\delta'}{s_{k+L-1}} \EPS(s_{k+L-1},s_{k+L})\EPS(t_{k+L-1},s_{k+L-1})\JS{s_{k+L-1}}{t_{k+L-1}},
\ED
where $\delta'=(-1)^{k+L-1}$. However, $L$ is even, so $\delta'=\delta$.
Further, all quantities on the right hand side of (\ref{2.6}) except the last one depend only on $\delta$ and the residue classes of $s_k$, $s_{k-1}$ and $t_{k-1}$ mod 4,
so we may write
\BD
\JS {t_{k+L}}{s_{k+L}} = \JS{-\delta}{s_{k}}\JS{\delta}{s_{k-1}} \EPS(s_{k-1},s_{k})\EPS(t_{k-1},s_{k-1})\JS{s_{k+L-1}}{t_{k+L-1}}.
\ED
Since $L$ is a period length of the Jacobi sequence, we see that the right hand side of this identity coincides with the right hand side of (\ref{2.6}).
Thus, $\JS {t_{k+L}}{s_{k+L}}=\JS {t_{k}}{s_{k}}$.

{\em Case 2}: $s$ is even.
\\ Since $t$ and $m$ are odd, both $s+t$ and $m+N$ are odd. The said theorem gives
\BE
\label{2.8}
\JS{-\delta s}{s+t}\JS pq\JS{\delta m}{m+N}=\EPS(s+t,q,m+N)
\EE
Here we use quadratic reciprocity and obtain
\BE
\label{2.10}
 \JS{m}{m+N}=\EPS(m,m+N)\JS{m+N}m=\EPS(m,m+N)\JS{N}m.
\EE
Similarly,
\begin{eqnarray}
\label{2.12}
  \JS s{s+t}=\JS{-t}{s+t}         &=&\JS{-1}{s+t}\JS t{s+t}=\\
  \JS{-1}{s+t}\EPS(t,s+t)\JS{s+t}t&=&\JS{-1}{s+t}\EPS(t,s+t)\JS st.\nonumber
\end{eqnarray}
From (\ref{2.8}), (\ref{2.10}) and (\ref{2.12}) we obtain an expression for $\JS Nm=\JS{t_k}{s_k}$ which depends only on the residue classes
of $s_k$, $t_k$, $s_{k-1}$ $t_{k-1} \mod 4$ and on $\JS st=\JS{s_{k-1}}{t_{k-1}}$. The same is true for $\JS{t_{k+L}}{s_{k+L}}$, the residue classes
of $s_ {k+L}$, $t_{k+L}$, $s_{k+L-1}$, $t_{k+L-1}\mod 4$ and $\JS{s_{k+L-1}}{t_{k+L-1}}$. Since $L$ is a period length of the Jacobi sequence,
we see that $\JS{t_{k+L}}{s_{k+L}}$ equals $\JS{t_k}{s_k}$.
\STOP

\section*{3. Proof of Theorem \ref{t1}}

As above, let $z=[\OV{a_0,a_1,\LD, a_{l-1}}]$ be a purely periodic quadratic irrational, the convergents $s_k/t_k$ being defined as in (\ref{1.0}) and (\ref{1.1}).
Let $L$ be an even multiple of $l$ such that $L$ is a period length of the Jacobi sequence of $z$ and (\ref{1.4}) holds. Again, we write
\BE
\label{3.2}
D_L=I+2^mU
\EE
with $m\ge 2$, $U$
as in (\ref{1.6}) such that not all entries of $U$ are even and $u=2^eu'$, $e\ge 0$, $u'$ odd. Let $k\ge 0$ be such that
$s_k\equiv 3 \mod 4$ and $t_k=2^{m+f}t'$ with $-m+1\le f\le e-1$, $t'$ odd. Note that $t_k$ is even but $s_k/t_k$ is not critical with respect to $L$ in the case $k\le L-1$.

\begin{lemma} 
\label{l1}
In the above setting, let $f\le e-2$. Then
\BD
  \JS{s_{k+L}}{t_{k+L}}=\JS{s_k}{t_k}
\ED
and $t_{k+L}=2^{m+f}t''$ with $t''\equiv t'\mod 4$.
In the case $f=e-1$ we have
\BD
 \JS{s_{k+L}}{t_{k+L}}=-\JS{s_k}{t_k}
\ED
and $t_{k+L}=2^{m-e-1}t''$ with $t''\equiv t'+2\mod 4$.

\end{lemma} 

\noindent
{\em Proof.}
From (\ref{2.0}), (\ref{3.2}) and  (\ref{1.6})  we obtain
\BE
\label{3.4}
 t_{k+L}=2^mus_k+t_k+2^mvt_k.
\EE
Since $u=2^eu'$ and $t_k=2^{m+f}t'$,
this reads
\BE
\label{3.6}
t_{k+L}=2^{m+f}(t'+2^{e-f}u's_k+ 2^mvt').
\EE
If $f\le e-2$, we obtain $t_{k+L}=2^{m+f}t''$ with  $t''\equiv t'\mod 4$ (observe $m\ge 2$).
Now the reciprocity law (\ref{1.2}) yields
\BD
\JS{s_{k+L}}{t_{k+L}}=\EPS(s_{k+L},t'')\JS{t_{k+L}}{s_{k+L}}.
\ED
The Kronecker symbol on the right hand side coincides with the Jacobi symbol, since $s_k$ and, consequently, $s_{k+L}$ is odd.
Moreover, $s_{k+L}\equiv s_k\mod 4$ and $t''\equiv t'$ mod 4. In addition, the reciprocal Jacobi sequence has the period length $L$.
From this we conclude
\BD
\JS{s_{k+L}}{t_{k+L}}=\EPS(s_{k},t')\JS{t_{k}}{s_{k}}.
\ED
On applying the reciprocity law (\ref{1.2}) again, we have
\BE
\label{3.8}
\JS{s_{k+L}}{t_{k+L}}=\JS{s_{k}}{t_{k}}.
\EE
In the case $f=e-1$ we observe that $e-f=1$ and $u's_k$ is odd. Accordingly, (\ref{3.6}) shows $t_{k+L}=2^{m+e-1}t''$ with $t''\equiv t'+2$ mod $4$.
Moreover, $s_{k+L}\equiv s_k\equiv 3\mod 4$, and so $\EPS(s_{k+L},t'')=-\EPS(s_k,t')$. This produces a sign change on the
right hand side of (\ref{3.8}). \STOP

\MN
The periodic case of the Kronecker symbol is contained in the following proposition.

\begin{prop} 
\label{p2}

In the above setting, suppose there are no critical convergents with respect to $L$. Then the Kronecker sequence is purely periodic
with period length $L$ except there is a convergent $s_k/t_k$ with $k\le L-1$, $s_k\equiv 3 \mod 4$ and $t_k=2^{m+e-1}t'$, $t'$ odd.
In this case the Kronecker sequence is purely periodic with period length $2L$.

\end{prop} 

\noindent
{\em Proof.} We consider an arbitrary convergent $s_k/t_k$. If $t_k$ is odd, the Kronecker symbol coincides with the Jacobi symbol, which means
$\JS{s_{k+L}}{t_{k+L}}=\JS{s_{k}}{t_{k}}$. If $t_k$ is even and $s_k\equiv 1\mod 4$, we have
\BD
 \JS{s_k}{t_k}=\EPS(s_k,t')\JS{t_k}{s_k},
\ED
where $t'$ is the odd part of $t_k$. However, $\EPS(s_k, t')=1$ since $s_k\equiv 1\MOD 4$. Hence we obtain $\JS{s_{k}}{t_{k}}=\JS{t_{k}}{s_{k}}$.
In the same way, $\JS{s_{k+L}}{t_{k+L}}=\JS{t_{k+L}}{s_{k+L}}$, because $s_{k+L}\equiv s_k \mod 4$. Now the periodicity of the reciprocal Jacobi sequence
$\JS{t_k}{s_k}$, $k\ge 0$,  shows $\JS{s_{k+L}}{t_{k+L}}=\JS{s_{k}}{t_{k}}$.

The case $t_k$ even and $s_k\equiv 3\mod 4$ is contained in Lemma 1. If $k\le L-1$ and $f\le e-2$, we have $\JS{s_{k+dL}}{t_{k+dL}}=\JS{s_{k}}{t_{k}}$
for all natural numbers $d$. Finally, if $k\le L-1$ and $f=e-1$, we obtain
\BD
 \JS{s_{k+dL}}{t_{k+dL}}=(-1)^d\JS{s_{k}}{t_{k}}.
\ED
In this situation the period length becomes $2L$.
\STOP

\MN
Suppose now that $k\le L-1$ and $s_k/t_k$ is critical with respect to $L$. Hence we have
$s_k\equiv 3\mod 4$ and $t_k\equiv 0\mod 2^{m+e}$. Recall the definition of $m$ and $e$: By (\ref{3.2}), we have $D_L=I+2^mU$, not all entries of $U$ even;
and the left lower entry $u$ of $U$ satisfies $u=2^eu'$, $u'$ odd. The relation (\ref{2.0}) implies $D_{2L}=D_L^2$. Therefore, this matrix reads
\BE
\label{3.7}
  D_{2L}=I+2^{m+1}\WT U \MB{ with } \WT U=U+2^{m-1}U^2.
\EE
Since $U$ is as in (\ref{1.6}), the left lower entry of $U^2$ equals $u(x+v)$, in particular, it is $\equiv 0 \mod 2^e$.
Accordingly, the left lower entry $\WT u$ of $\WT U$, i.e., $\WT u=u+2^{m-1}u(x+v)$, is $\equiv u\mod 2^{e+1}$, since $m\ge 2$.
Hence $\WT u=2^eu''$, $u''$ odd, with the same exponent $e$.
This means the following. If $t_k\equiv 0 \mod 2^{m+e+1}$, the convergent
$s_k/t_k$ is also critical with respect to $2L$.

Hence there is a number $r\ge 0$ such that $s_k/t_k$ is critical with
respect to $2^rL$ but not with respect to $2^{r+1}L$. In the following lemma we suppose that $r$ has been chosen in this way.
For the sake of simplicity, however, we simply write $L$ instead
of $2^rL$ and adopt the other notations connected with $D_L$. Then $t_k=2^{m+e}t'$, $t'$ odd.

\begin{lemma} 
\label{l2}
In the above setting, let $k\le L-1$. Suppose that $s_k/t_k$ is critical with respect to $L$ and $t_k=2^{m+e}t'$, $t'$ odd.
Suppose, further, that $s_k/t_k$ is not critical with respect to $2L$.
Then $s_{k+L}/t_{k+L}$ is critical with respect to $2L$. Moreover, for every $d\ge 1$,
\BD
 \JS{s_{k+d\cdot 2L}}{t_{k+d\cdot 2L}}=(-1)^d\JS{s_{k}}{t_{k}}.
\ED
\end{lemma}

\noindent
{\em Proof.} Our situation corresponds to the case $f=e$ in formula (\ref{3.6}), and so
\BD
 t_{k+L}=2^{m+e}(t'+u's_k+ 2^mvt').
\ED
Here, however $t'$, $u'$ and $s_k$ are odd. Accordingly, $t_{k+L}\equiv 0 \mod 2^{m+e+1}$, which means that $s_{k+L}/t_{k+L}$ is critical with respect to $2L$,
as we have seen above.

As in (\ref{3.7}) we write $D_{2L}=D_L^2=I+2^{m+1}\WT U$, where the matrix $\WT U$ has the lower entries $\WT u=2^eu''$, $u''$ odd, and $\WT v$.
Here the analogue of (\ref{3.4}) reads
\BE
\label{3.9}
  t_{k+2L}=2^{m+1}\WT us_k+t_k+2^{m+1}\WT vt_k.
\EE
If we insert $\WT u=2^eu''$ and $t_k=2^{m+e}t'$, we obtain
\BD
  t_{k+2L}=2^{m+e}(2u''s_k+t'+2^{m+1}\WT vt').
\ED
Because $u''$ and $s_k$ are odd, this yields
\BE
\label{3.10}
 t_{k+2L}=2^{m+e}t''\MB{ with }t''\equiv t'+2\mod 4.
\EE
As in the proof of Lemma \ref{l1} we use the reciprocity law and obtain
\BD
  \JS{s_{k+2L}}{t_{k+2L}}=\EPS(s_{k+2L},t'')\JS{t_{k+2L}}{s_{k+2L}}.
\ED
By (\ref{3.10}), $\EPS(s_{k+2L},t'')=-\EPS(s_k,t')$. Now the periodicity of the reciprocal Jacobi sequence combined with another application of the reciprocity law
yields
\BD
  \JS{s_{k+2L}}{t_{k+2L}}=-\JS{s_{k}}{t_{k}}.
\ED
Next formula (\ref{3.9}) gives
\BD
  t_{k+4L}=2^{m+1}\WT us_{k+2L}+t_{k+2L}+2^{m+1}\WT vt_{k+2L}.
\ED
On inserting $\WT u=2^eu''$ and $t_{k+2L}=2^{m+e}t''$, we have
\BD
t_{k+4L}=2^{m+e}(2u''s_{k+2L}+t''+2^{m+1}\WT vt'').
\ED
Since $u''$ and $s_{k+2L}$ are odd, this gives
\BD
 t_{k+4L}=2^{m+e}t'''\MB{ with }t'''\equiv t''+2\mod 4.
\ED
Now the above arguments show
\BD
  \JS{s_{k+4L}}{t_{k+4L}}=-\JS{s_{k+2L}}{t_{k+2L}}=\JS{s_k}{t_k}.
\ED
The general case $k+d\cdot 2L$, $d\ge 1$, is settled in the same way by induction.
\STOP

\MN
Let $L$ be such that $s_k/t_k$ is critical with respect to $L$. Then there is a number $r_1\ge 0$ such that $s_k/t_k$ is critical with respect to $2^{r_1}L$
but not with respect to $2^{r_1+1}L$.
Put $k_1=k$ and $k_2=k+2^{r_1}L$. By Lemma \ref{l2}, $s_{k_2}$ is critical with respect to $2^{r_1+1}L$. Hence there is a number $r_2>r_1$ such that
$s_{k_2}/t_{k_2}$ is critical with respect to $2^{r_2}L$ but not with respect to $2^{r_2+1}L$. In this way we obtain an infinite sequence $(k_j,r_j)$, $j\ge 1$, which is strictly increasing in both arguments such
that $s_{k_j}/t_{k_j}$ is critical with respect to $2^{r_j}L$ but not with respect to $2^{r_j+1}L$.

\MN
{\em Example.} In the case of the above example $(7+\sqrt{82})/11=[\OV{1,2,5}]$ the convergent $s_7/t_7$ is critical with respect to $L=12$, but not critical with respect to $2L=24$.
Hence we have $k_1=7$, $r_1=0$. Since $7+L=19$,
Lemma \ref{l2} says that $s_{19}/t_{19}$ is critical with respect to $2L$. Since it is not critical with respect to $4L$, we have $k_2=19$, $r_2=1$. Now $19+2L=43$, and
$s_{43}/t_{43}$ is critical with respect to $8L$, but not with respect to $16L$. Hence $k_3=43$ and $r_3=3$. Next we have $43+8L=139$,
$s_{139}/t_{139}$ is critical with respect to $64L$, but not with respect to $128L$. So we have $k_4=139$ and $r_4=6$. Accordingly, the first members
of our sequence are $(7,0), (19,1), (43,3)$ and $(139, 6)$.

\MN
{\em Proof of Theorem \ref{t1}.} We have only to consider the case that there is a critical convergent with respect to $L$. Hence we know that there is a sequence
$(k_j,r_j)$ with the above properties. Suppose that the Kronecker sequence is periodic with period length $2^idL$, $i\ge 0$, $d\ge 1$, $d$ odd.
Then there is an integer $k_0\ge 0$ such that for all $k\ge k_0$
\BD
\JS{s_{k+2^idL}}{t_{k+2^idL}}=\JS{s_k}{t_k}.
\ED
We choose a member $(k_j,r_j)$ of our sequence in such a way that $k_j\ge k_0$ and $r_j\ge i-1$. Then $2^{r_j+1}dL$ is a multiple of $2^idL$,
and so
\BD
\JS{s_{k_j+2^{r_j+1}dL}}{t_{k_j+2^{r_j+1}dL}}=\JS{s_{k_j}}{t_{k_j}}.
\ED
By Lemma \ref{l2}, however,
\BD
 \JS{s_{k_j+2^{r_j+1}dL}}{t_{k_j+2^{r_j+1}dL}}=(-1)^d\JS{s_{k_j}}{t_{k_j}}.
\ED
Since $d$ is odd, this is a contradiction.\STOP


\vspace{0.5cm}
\noindent
Kurt Girstmair            \\
Institut f\"ur Mathematik \\
Universit\"at Innsbruck   \\
Technikerstr. 13/7        \\
A-6020 Innsbruck, Austria \\
Kurt.Girstmair@uibk.ac.at

\end{document}